\documentclass[11pt]{amsart}
\usepackage{amsmath,amssymb,latexsym,bm,xcolor}

\newcommand{\mcr}[1]{\lfloor #1 \rfloor}

\newcommand{\bZ}{\mbox{${\mathbb Z}$}}

\newcommand{\bC}{\mbox{${\mathbb C}$}}

\newcommand{\QQ}{\mbox{${\mathbf Q}_G$}}
\newcommand{\OQGl}[1]{[{\mathcal O}_{{\mathbf Q}_G}(#1)]}
\newcommand{\OQG}[1]{[{\mathcal O}_{{\mathbf Q}_G(#1)}]}

\def\mytilde{\kern-.015in\hbox{\lower.03in\hbox{\~{}}}\kern-.01in}
\def\Waf{W_{\mathrm{af}}}
\def\Wafp{W_{\mathrm{af}}^{\ge 0}}

\DeclareMathOperator{\wend}{end}
\DeclareMathOperator{\hgt}{height}
\DeclareMathOperator{\wt}{wt}
\DeclareMathOperator{\dn}{down}

\newcommand{\QLS}{\mathop{\rm QLS}\nolimits}
\newcommand{\QB}{\mathop{\rm QB}\nolimits}
\newcommand{\Deg}{\mathop{\rm deg}\nolimits}

\newcommand{\rh}{\widehat{r}}

\newtheorem{thm}{Theorem}

\newtheorem{prop}[thm]{Proposition}

\newtheorem{dfn}[thm]{Definition}

\newtheorem{exa}[thm]{Example}

\newtheorem{rema}[thm]{Remark}

\makeatletter

\let\choose\@@choose

\makeatother

\usepackage[left=2.54cm, right=2.54cm, top=2.54cm, bottom=2.54cm]{geometry}

\title[A Chevalley formula for semi-infinite flag manifolds]{A Chevalley formula for semi-infinite\\ flag manifolds and quantum $K$-theory\\  {\rm{\small{Extended Abstract}}}}

\author[C.~Lenart]{Cristian Lenart}
\address[Cristian Lenart]{Department of Mathematics and Statistics, State University of New York at Albany, 
Albany, NY 12222, U.S.A.}
\email{clenart@albany.edu}

\author[S.~Naito]{Satoshi Naito}
\address[Satoshi Naito]{Department of Mathematics, Tokyo Institute of Technology,
2-12-1 Oh-Okayama, Meguro-ku, Tokyo 152-8551, Japan}
\email{naito@math.titech.ac.jp}

\author[D.~Sagaki]{Daisuke Sagaki}
\address[Daisuke Sagaki]{Institute of Mathematics, University of Tsukuba, 
Tsukuba, Ibaraki 305-8571, Japan}
\email{sagaki@math.tsukuba.ac.jp}


\begin{document}

\begin{abstract} We give a combinatorial Chevalley formula for an arbitrary weight, in the torus-equivariant $K$-theory of semi-infinite flag manifolds, which is expressed in terms of the quantum alcove model. As an application, we prove the Chevalley formula for anti-dominant fundamental weights in the (small) torus-equivariant quantum $K$-theory $QK_{T}(G/B)$ of the flag manifold $G/B$; this has been a longstanding conjecture about the multiplicative structure of $QK_{T}(G/B)$. Moreover, in type $A_{r}$, we prove that the so-called quantum Grothendieck polynomials indeed represent Schubert classes in the (non-equivariant) quantum $K$-theory $QK(SL_{r+1}/B)$.
\end{abstract}

\maketitle

\section{Introduction}

This paper is concerned with a geometric application of the combinatorial model known as the {\em quantum alcove model}, introduced in \cite{lalgam}. Its precursor, the alcove model of the first author and Postnikov, was used to uniformly describe the  highest weight {\em Kashiwara crystals} of symmetrizable Kac-Moody algebras \cite{LP1}, as well as the {\em Chevalley formula} in the equivariant $K$-theory of flag manifolds $G/B$ \cite{LP}. More generally, the quantum alcove model was used to uniformly describe certain crystals of affine Lie algebras (single-column {\em Kirillov-Reshetikhin crystals}) and {\em Macdonald polynomials} specialized at $t=0$ \cite{lnsumk2,lnsumk3}. The objects of the quantum alcove model (indexing the crystal vertices and the terms of Macdonald polynomials) are paths in the {\em quantum Bruhat graph} on the Weyl group~\cite{BFP}. In this paper we complete the above picture, by extending to the quantum alcove model the geometric application of the alcove model, namely the $K$-theory Chevalley formula. 

To achieve our goal, we need to consider the so-called {\em semi-infinite flag manifold} $\QQ$. 
We give a Chevalley formula for an arbitrary weight in the {$T \times \mathbb{C}^{*}$}-equivariant $K$-group $K_{T\times{\mathbb C}^*}(\QQ)$ of $\QQ$, which is described in terms of the quantum alcove model.
In~\cite{knsekt} and \cite{nospcf}, the Chevalley formulas for $K_{T\times{\mathbb C}^*}(\QQ)$ were originally given in terms of the {\em quantum LS path model} in the case of a dominant and an anti-dominant weight, respectively.
For a general (not dominant nor anti-dominant) weight, there is no quantum LS path model, but there is a quantum alcove model.
Hence, in order to obtain a Chevalley formula for an arbitrary weight, we first need to translate the formulas above to the quantum alcove model by using the weight-preserving bijection between the two models given by Proposition~\ref{qam2qls}.
Based on these translated formulas (Theorems~\ref{chevdomqam} and \ref{corantidom}), we obtain a Chevalley formula (Theorem~\ref{genchev}) for an arbitrary weight by using elaborate combinatorics of the quantum alcove model.

The study of the equivariant $K$-group of semi-infinite flag manifolds was started in~\cite{knsekt}, and a breakthrough in this study is \cite{kat1} (see also \cite{kat2}), in which Kato established a $\mathbb{C}[P]$-module isomorphism from the (small) $T$-equivariant quantum $K$-group $QK_{T}(G/B)$ of the finite-dimensional flag manifold $G/B$ onto (a version of) the $T$-equivariant $K$-group $K_{T}^{\prime}(\QQ)$ of $\QQ$, where $P$ is the weight lattice generated by the fundamental weights $\varpi_{i}$, $i \in I$.
This isomorphism sends each (opposite) Schubert class in $QK_{T}(G/B)$ to the corresponding semi-infinite Schubert class in $K_{T}^{\prime}(\QQ)$; moreover, it respects the quantum multiplication in $QK_{T}(G/B)$ and the tensor product in $K_{T}^{\prime}(\QQ)$. Based on this result, a longstanding conjecture on the multiplicative structure of $QK_{T}(G/B)$, i.e., the Chevalley formula (Theorem~\ref{qkchev}) for anti-dominant fundamental weights $- \varpi_{i}$, $i \in I$, for $QK_{T}(G/B)$ {is proved by our Chevalley formula for $K_{T \times \mathbb{C}^{*}}(\QQ)$ under the specialization at $q = 1$.} 

As another application of our Chevalley formula, we can prove an important conjecture for the non-equivariant quantum $K$-group $QK(SL_{r+1}/B)$ of the flag manifold of type $A_{r}$ (Theorem~\ref{qgroth}):
the {\em quantum Grothendieck polynomials}, introduced in \cite{lamqgp}, indeed represent Schubert classes in $QK(SL_{r+1}/B)$. In this way, we generalize  the results of~\cite{fgpqsp}, where the {\em quantum Schubert polynomials} are constructed as representatives for Schubert classes in the quantum cohomology of $SL_{r+1}/B$.  
Therefore, we can use quantum Grothendieck polynomials to compute any structure constant in $QK(SL_{r+1}/B)$ (with respect to the Schubert basis);  
indeed, we just need to expand their products in the basis they form, which is done by \cite[Algorithm~3.28]{lamqgp}, see \cite[Example~7.4]{lamqgp}. This is important, since computing even simple products in quantum $K$-theory is notoriously difficult.

\subsection*{Acknowledgments}

 C.L. was partially supported by the NSF grants DMS-1362627 and DMS-1855592.
 S.N. was partially supported by JSPS Grant-in-Aid for Scientific Research (B) 16H03920.
 D.S. was partially supported by JSPS Grant-in-Aid for Scientific Research (C) 15K04803 and 19K03415.

\section{Background on the combinatorial models}\label{backg}

\subsection{Root systems} Let $\mathfrak{g}$ be a complex simple Lie algebra, and $\mathfrak{h}$ a Cartan subalgebra.  Let $\Phi\subset \mathfrak{h}^*$ be the corresponding irreducible \textit{root system}, $\mathfrak{h}^*_\mathbb{R}$ the real span of the roots, and $\Phi^{+}\subset\Phi$ the set of positive roots. 
As usual, we denote $\rho:=\frac{1}{2}(\sum_{\alpha\in\Phi^+}\alpha)$.
  Let $\alpha_i\in\Phi^+$ be the \textit{simple roots}, for $i$ in an indexing set $I$. We denote $\langle\cdot,\cdot\rangle$ the nondegenerate scalar product on $\mathfrak{h}^*_{\mathbb{R}}$ induced by the Killing form.  Given a root $\alpha$, we consider the corresponding \textit{coroot} $\alpha^{\vee}$ and reflection $s_{\alpha}$. The root and coroot lattices are denoted by $Q$ and $Q^\vee$, as usual, while the positive part of the coroot lattice is denoted by $Q^{\vee,+}$. The \textit{weight lattice} $P$ is generated by the \textit{fundamental weights} $\varpi_i$, for $i\in I$, which form the dual basis to the simple coroots. Let $P^+$ be the set of \textit{dominant weights}. 

Let $W$ be the \textit{Weyl group}, with length function $\ell(\cdot)$ and longest element $w_\circ$.  The \textit{Bruhat order} on $W$ is defined by its covers $w\lessdot ws_{\alpha},$ for $\ell(ws_\alpha) = \ell(w)+1$, where $\alpha\in\Phi^+$. 

Given $\alpha\in\Phi$ and $k\in\mathbb{Z}$, we denote by $s_{\alpha,k}$ the reflection in the affine hyperplane $H_{\alpha,k}:=\{\lambda\in\mathfrak{h}^*_{\mathbb{R}} \mid \langle\lambda,\alpha^{\vee}\rangle = k\}$. 
These reflections generate the \textit{affine Weyl group} $\Waf=W\ltimes Q^\vee$ for the \textit{dual root system} $\Phi^{\vee}$. The hyperplanes $H_{\alpha,k}$ divide the vector space $\mathfrak{h}^*_{\mathbb{R}}$ into open regions, called \textit{alcoves}.  The \textit{fundamental alcove} is denoted by $A_{\circ}$. 

 The \textit{quantum Bruhat graph} ${\QB}(W)$ on $W$ is defined by adding downward (quantum) edges, denoted $w \triangleleft ws_\alpha$, to the covers of the Bruhat order, i.e., the edges of ${\QB}(W)$ are:
$$w\xrightarrow{\alpha}ws_\alpha\hspace{8pt}\text{if}\hspace{6pt} w\lessdot ws_\alpha\hspace{6pt}\mbox{or}\hspace{6pt} \ell(ws_{\alpha}) =\ell(w) - 2\langle\rho,\alpha^{\vee}\rangle + 1\,,\hspace{8pt}\mbox{where $\alpha\in\Phi^+$}\,.$$\label{quantum bruhat graph eq}
We define the {\em weight} of an edge $w \overset{\alpha}\longrightarrow {ws_\alpha}$ to be
either $\alpha^\vee$ or $0$, depending on whether it is a quantum edge or not, respectively. 
Then the weight of a directed path is 
the sum of the weights of its edges. It turns out that the weight of a shortest directed path from $v$ to $w$ is independent of the mentioned path, so we will denote it by $\wt(w\Rightarrow v)$; see~\cite{lnsumk2}. 

For the remainder of this section, we fix $\lambda\in P^+$. Let $W_J$ be the stabilizer of $\lambda$, as a parabolic subgroup with $J\subset I$ and root system $\Phi_J$. We denote the set of minimum-length coset representatives for $W/W_J$ by $W^J$, and the minimum-length coset representative of $wW_J$ by $\mcr{w}$. We consider the {\em parabolic quantum Bruhat graph} on $W^J$, denoted by $\QB(W^J)$; this generalizes $\QB(W)$, see~\cite{lnsumk1}. Its directed edges are labeled by $\alpha\in \Phi^+\setminus \Phi_J^+$. The upward edges are the covers of the Bruhat order on $W^J$, while the downward (quantum) edges $w \overset{\alpha}{\longrightarrow} \mcr{ws_\alpha}$ are given by the condition $\ell(\mcr{ws_\alpha})=\ell(w)-2\langle \rho-\rho_J, \alpha^\vee \rangle+1$. 
Given a rational number $b$, we define $\QB_{b\lambda}(W^J)$ to be the subgraph of 
 $\QB(W^{J})$ with 
the same vertex set but having only the edges with labels $\alpha$ satisfying $b\langle\lambda,\alpha^\vee\rangle\in{\mathbb Z}$.

We now recall the quantum Bruhat graph analogue of a certain lift from $W/W_J$ to $W$ which was previously defined by Deodhar. Let $\ell(w \Rightarrow x)$ denote the length of the shortest path from $w$ to $x$ in $\QB(W)$. It was shown in \cite{lnsumk1} that, given $v,w\in W$, there exists a unique element $x \in vW_{J}$ such that 
$\ell(w \Rightarrow x)$ attains its minimum value as a function of $x\in vW_J$. For reasons explained in \cite{lnsumk1}, we denote the unique element by $\min(vW_J,\preceq_w)$, and call it a {\em quantum Deodhar lift}.

\subsection{Quantum LS paths}\label{sec:qls}

\begin{dfn}[\cite{lnsumk2}]\label{defls} A quantum LS path $\eta\in\QLS(\lambda)$, for $\lambda\in P^+$, is given by two sequences
\begin{equation}\label{E:lsste}
(0 = b_1<b_2<b_3<\cdots<b_t<b_{t+1}=1)\,;\;\; \qquad(\phi(\eta)=\sigma_1,\,\sigma_{2},\,\ldots,\,\sigma_t=\iota(\eta))\,,
\end{equation}
where $b_k\in \mathbb{Q}$, $\sigma_k\in W^J$, and there is a directed path in $\QB_{b_k\lambda}(W^J)$ from $\sigma_{k-1}$ to $\sigma_k$, for each $k=2,\ldots,t$. The elements $\sigma_k$ are called the {\em directions} of $\eta$, while $\iota(\eta)$ and $\phi(\eta)$ are the {\em initial} and {\em final directions}, respectively. 
\end{dfn}

This data encodes the sequence of vectors 
\begin{equation}\label{seqvec}v_t:=(b_{t+1}-b_t)\sigma_t\lambda\,, \;\:\ldots\,,\;\: v_2:=(b_3-b_2)\sigma_2\lambda\,,\;\: v_1:=(b_2-b_1)\sigma_1\lambda\,.\end{equation}
 We can view $\eta\in\QLS(\lambda)$ as a piecewise-linear 
path given by the sequence of points 
\[0\,,\;\: v_t\,,\;\: v_{t-1}+v_t\,,\;\: \dotsc\,,\;\: v_1+\cdots+v_t\,.\]
The endpoint of the path, also called its weight, is $\wt(\eta):=\eta(1)=v_1+\cdots+v_t$. 

Given $w\in W$, we define the {\em initial direction of $\eta$ with respect to $w$} as $\iota({\eta},{w}):=w_t \in W$, where the sequence $(w_k)$ is calculated by the following recursive formula: 
\begin{equation} \label{eq:tiw}
w_0:=w\,,\;\;\;\; w_k:=\min(\sigma_kW_J,\preceq_{w_{k-1}})\;\; \text{for $k=1,\ldots,t$}\,. 
\end{equation}
Also, we set
\begin{align} \label{eq:xiv}
\xi(\eta,w)&:= \sum_{k=1}^{t} \wt ({w}_{k-1} \Rightarrow {w}_{k})\,, \\
 \label{eq:degx}
\Deg_{w}(\eta)&:= - \sum_{k=1}^{t} (1-b_k) \langle{\lambda},{\wt({w}_{k-1} \Rightarrow {w}_{k})}\rangle\,.
\end{align}

\subsection{The quantum alcove model}\label{sec:qam}

We say that two alcoves are {adjacent} if they are distinct and have a common wall. Given a pair
of adjacent alcoves $A$ and $B$, we write $A \xrightarrow{\beta}  B$  for $\beta\in\Phi$ if the
common wall is orthogonal to $\beta$ and $\beta$ points in the direction from $A$ to $B$. 

\begin{dfn}[\cite{LP}]
	An  \emph{alcove path} is a sequence of alcoves $(A_0, A_1, \ldots, A_m)$ such that
	$A_{j-1}$ and $A_j$ are adjacent, for $j=1,\ldots, m.$ We say that $(A_0, A_1, \ldots, A_m)$  
	is \emph{reduced} if it has minimal length among all alcove paths from $A_0$ to $A_m$.
\end{dfn}

Let $\lambda\in P$ be any weight, although dominant and anti-dominant $\lambda$  will play a special role.
Let $A_{\lambda}=A_{\circ}+\lambda$ be the translation of the fundamental alcove $A_{\circ}$ by $\lambda$.
	
\begin{dfn}[\cite{LP}]\label{deflch}
	The sequence of roots $\Gamma(\lambda)=(\beta_1, \beta_2, \dots, \beta_m)$ is called a
	\emph{$\lambda$-chain} if 
	\[	
		A_0=A_{\circ} \xrightarrow{-\beta_1}  A_1
		\xrightarrow{-\beta_2} \cdots 
		\xrightarrow{-\beta_m}  A_m=A_{-\lambda}
	\]
is a reduced alcove path.
\end{dfn}

A reduced alcove path $(A_0=A_{\circ},A_1,\ldots,A_m=A_{-\lambda})$ can be identified 
with the corresponding total order on the hyperplanes, to be called 
{\em $\lambda$-hyperplanes}, which separate $A_\circ$ from $A_{-\lambda}$. This total order is given by the sequence 
$H_{\beta_i,-l_i}$ for $i=1,\ldots,m$, where $H_{\beta_i,-l_i}$ contains the common wall of 
$A_{i-1}$ and $A_i$. Note that $\langle\lambda,\beta_i^\vee\rangle\ge0$, and that the integers $l_i$, called {\em heights}, have the following ranges:
\begin{equation}\label{ranges}0\le l_i\le\langle\lambda,\beta_i^\vee\rangle-1\;\;\mbox{if}\;\;\beta_i\in\Phi^+\,, \;\;\;\;\mbox{and}\;\;\;\;1\le l_i\le\langle\lambda,\beta_i^\vee\rangle\;\;\mbox{if}\;\;\beta_i\in\Phi^-\,.\end{equation}
 Note also that a $\lambda$-chain $(\beta_1, \ldots, \beta_m)$ determines 
the corresponding reduced alcove path, so we can identify them as well. 

\begin{rema}\label{remredword} {\rm 
An alcove path corresponds to the choice of a reduced word for the affine Weyl group element sending $A_\circ$ to $A_{-\lambda}$ \cite[Lemma 5.3]{LP}. 
}
\end{rema}

For dominant $\lambda$, we have the following special choice of a $\lambda$-chain, by~\cite[Section 4]{LP1}. 

\begin{prop}[\cite{LP1}]
\label{speciallch} 
Given a total order $I=\{1<2<\dotsm<r\}$ on the Dynkin nodes, we express $\beta^\vee=\sum_{i=1}^r c_i \alpha_i^\vee$
in the ${\mathbb Z}$-basis of simple coroots. Consider the
total order on the set of $\lambda$-hyperplanes defined by 
the lexicographic order on their images in ${\mathbb Q}^{r+1}$ under the map
\begin{equation}
\label{E:stdvec}
	H_{\beta,-l}\mapsto \frac{1}{\langle\lambda, \beta^\vee \rangle} (l,c_1,\ldots,c_r).
\end{equation}
This map is injective, thereby endowing 
the set of $\lambda$-hyperplanes with a total order, which is a $\lambda$-chain. We call it the
{\em lexicographic (lex) $\lambda$-chain}, and denote it by $\Gamma_{\rm lex}(\lambda)$. 
\end{prop}

The rational number $l/\langle\lambda, \beta^\vee \rangle$ is called the {\em relative height} of the $\lambda$-hyperplane $H_{\beta,-l}$. By definition, the sequence of relative heights in the lex $\lambda$-chain is weakly increasing. 

The objects of the quantum alcove model are defined next; for examples, we refer to~\cite{lalgam,lnsumk2}. Compared to the original construction in~\cite{lalgam}, here we consider a generalization of this model, by letting $\lambda$ be any weight, as opposed to only a dominant weight; another aspect of the generalization is making the model depend on a fixed element  $w\in W$, such that the initial model corresponds to $w$ being the identity. In addition to $w$, we fix an arbitrary $\lambda$-chain $\Gamma(\lambda)=(\beta_1,\,\ldots,\,\beta_m)$, and let $r_i:=s_{\beta_i}$, $\widehat{r}_i:=s_{\beta_i,-l_i}$. 

\begin{dfn}[\cite{lalgam}]
\label{def:admissible}
	A subset 
	$A=\left\{ j_1 < \cdots < j_s \right\}$ of $[m]:=\{1,\ldots,m\}$ (possibly empty)
 	is a $w$-\emph{admissible subset} if
	we have the following directed path in $\QB(W)$:
	\begin{equation}
	\label{eqn:admissible}
	 \Pi(w,A):\;\;\;\;w\xrightarrow{|\beta_{j_1}|} w r_{{j_1}} 
	\xrightarrow{|\beta_{j_2}|}  wr_{{j_1}}r_{{j_2}} 
	\xrightarrow{|\beta_{j_3}|}  \cdots 
	\xrightarrow{|\beta_{j_s}|}  wr_{{j_1}}r_{{j_2}} \cdots r_{{j_s}}=:\wend(w,A)\,.
	\end{equation}
 	We let ${\mathcal A}(w,\Gamma(\lambda))$ be the collection of all $w$-admissible subsets of $[m]$.
\end{dfn}

We now associate several parameters with the pair $(w,A)$. 
The weight of $(w,A)$ is  
	\begin{equation}
	\label{defwta}
	\wt(w,A):=-w\rh_{{j_1}} \cdots 
	         \rh_{{j_s}}(-\lambda)\,.
	\end{equation}

Given the height sequence $(l_1,\ldots,l_m)$ above, we define 
the complementary height sequence $(\widetilde{l}_1,\ldots,\widetilde{l}_m)$ 
by $\widetilde{l}_i:=\langle\lambda,\beta_i^\vee\rangle-l_i$. 
Given $A=\{j_1<\cdots<j_s\}\in{\mathcal A}(w,\Gamma(\lambda))$, let 
\begin{equation*}
A^{-}:=\bigl\{j_i \in A \mid 
  wr_{{j_1}} \cdots r_{{j_{i-1}}} > 
 w r_{{j_1}} \cdots r_{{j_{i-1}}}r_{{j_{i}}}
\bigr\}\,;
\end{equation*}
in other words, we record the quantum steps 
in the path $\Pi(w,A)$ defined in \eqref{eqn:admissible}. 
Let
\begin{equation}
\label{defheight}
\dn(w,A):=\sum_{j\in A^-}|\beta_j|^\vee\in Q^{\vee,+}\,,\;\;\;\;\hgt(w,A):=\sum_{j\in A^-}\widetilde{l}_j\,.
\end{equation}

\section{Chevalley formulas for the semi-infinite flag manifold}\label{csi}

Consider a simply-connected simple algebraic group $G$ over $\mathbb{C}$, with Borel subgroup $B = T N$, maximal torus $T$, and unipotent radical $N$. The full {\em semi-infinite flag manifold} $\mathbf{Q}_{G}^{\mathrm{rat}}$ is the reduced (ind-)scheme associated to $G(\bC((z))\,)/(T\cdot N(\bC((z)))\,)$; in this paper, we concentrate on its semi-infinite Schubert subvariety $\QQ := \QQ(e) \subset \mathbf{Q}_{G}^{\mathrm{rat}}$ corresponding to the identity element $e \in \Waf$, which we also call the semi-infinite flag manifold. The $T \times \mathbb{C}^*$-equivariant $K$-group $K_{T \times \mathbb{C}^*}(\QQ)$ of $\QQ$ 
{has a (topological) $\mathbb{C}[q, q^{-1}][P]$-basis} of {\em semi-infinite Schubert classes}, and its multiplicative structure is determined by a {\em Chevalley formula}, which expresses the tensor product of a Schubert class with the class of a line bundle. In \cite{knsekt} and \cite{nospcf}, the Chevalley formulas were given in the case of a dominant and an anti-dominant weight $\lambda$, respectively. These formulas were expressed in terms of the quantum LS path model. We will express them in terms of the quantum alcove model based on the lexicographic $\lambda$-chain. The goal is to generalize these formulas for an arbitrary weight $\lambda$, and we will also see that an arbitrary $\lambda$-chain can be used. Throughout this section, $W_J$ is the stabilizer of $\lambda$, and we use freely the notation in Section~\ref{backg}.

The $T \times \mathbb{C}^*$-equivariant $K$-group $K_{T \times \mathbb{C}^*}(\QQ)$ is a $\mathbb{C}[q, q^{-1}][P]$-module consisting of all (possibly infinite) linear combinations of the classes $\OQG{x}$ of the structure sheaves of the semi-infinite Schubert varieties $\QQ(x) (\subset \QQ)$ with coefficients $a_{x} \in \mathbb{C}[q, q^{-1}][P]$ for $x \in \Wafp = W \times Q^{\vee,+}$ such that the sum $\sum_{x \in \Wafp} \vert a_{x} \vert$ of the absolute values $\vert a_{x} \vert$ lies in $\mathbb{C}(\!(q^{-1})\!)[P]$. Here $\mathbb{C}^*$ acts on $\QQ$ by loop rotation, and $\mathbb{C}[P]$ ($= \mathbb{C} \otimes_{\mathbb{Z}} \mathbb{Z}[P]$) is the group algebra of $P$, spanned by formal exponentials $\mathbf{e}^\lambda$, for $\lambda\in P$, with $\mathbf{e}^\lambda\mathbf{e}^\mu=\mathbf{e}^{\lambda+\mu}$; note that $\mathbb{Z}[P]$ is identified with the representation ring of $T$.
We also consider the $\mathbb{C}[q, q^{-1}][P]$-submodule $K_{T \times \mathbb{C}^*}^{\prime}(\QQ)$ of $K_{T \times \mathbb{C}^*}(\QQ)$ consisting of all finite linear combinations of the classes $\OQG{x}$ with coefficients in $\mathbb{C}[q, q^{-1}][P]$ for $x \in \Wafp$.
The $T$-equivariant $K$-groups of $\QQ$, denoted by $K_{T}(\QQ)$ and $K_{T}^{\prime}(\QQ)$, are obtained from the $K_{T \times \mathbb{C}^*}(\QQ)$ and $K_{T \times \mathbb{C}^*}^{\prime}(\QQ)$ above, respectively, by the specialization $q = 1$. Hence the Chevalley formulas in $K_T(\QQ)$ (for anti-dominant weights) and $K_{T}^{\prime}(\QQ)$ (for arbitrary weights) are obtained from the corresponding one in $K_{T \times \mathbb{C}^*}(\QQ)$ by setting $q=1$.

\subsection{Chevalley formulas for dominant and anti-dominant weights} 
We start with the Chevalley formula for dominant weights, which was derived in terms of semi-infinite LS paths in~\cite{knsekt}, and then restated in \cite[Corollary~C.3]{nospcf} in terms of quantum LS paths. 

Let $\lambda = \sum_{i \in I} \lambda_i \varpi_i$ be a dominant weight. We denote by $\overline{{\rm Par}(\lambda)}$ the set of $I$-tuples of partitions ${\bm{\chi}} = (\chi^{(i)})_{i \in I}$ such that $\chi^{(i)}$ is a partition of length at most $\lambda_i$ for all $i \in I$. For ${\bm{\chi}} = (\chi^{(i)})_{i \in I} \in \overline{{\rm Par}(\lambda)}$, we set $|{\bm{\chi}}| := \sum_{i \in I} |\chi^{(i)}|$, with $|\chi^{(i)}|$ the size of the partition $\chi^{(i)}$. Also set $\iota({\bm{\chi}}) := \sum_{i \in I} \chi^{(i)}_1 \alpha_i^{\vee} \in Q^{\vee,+}$, with $\chi^{(i)}_1$ the first part of the partition $\chi^{(i)}$.

\begin{thm}[\cite{knsekt,nospcf}]\label{chevdomqls} Let $x=wt_{\xi}\in \Wafp = W \times Q^{\vee,+}$. Then, in $K_{T \times \mathbb{C}^*}(\QQ)$, we have
\[\begin{split}&\OQGl{-w_\circ\lambda}\cdot\OQG{x}=\\[3mm]&\qquad=\sum_{\eta\in\QLS(\lambda)}\sum_{{\bm{\chi}}\in\overline{{\rm Par}(\lambda)}}q^{\Deg_w(\eta)-\langle\lambda,\xi\rangle-|{\bm{\chi}}|}\mathbf{e}^{\wt(\eta)}\OQG{\iota(\eta,w)t_{\xi+\xi(\eta,w)+\iota({\bm{\chi}})}}\,.\end{split}\]
\end{thm}

We now translate this formula in terms of the quantum alcove model for the lex $\lambda$-chain $\Gamma_{\rm lex}(\lambda)$. To this end, given $w\in W$, we construct a bijection $A\mapsto \eta$ between ${\mathcal A}(w,\Gamma_{\rm lex}(\lambda))$ and $\QLS(\lambda)$, for which several properties are then proved. 

In order to construct the forward map, let $A=\{j_1<\cdots<j_s\}$ be in ${\mathcal A}(w,\Gamma_{\rm lex}(\lambda))$. The corresponding heights are within the first range in~\eqref{ranges}. Consider the weakly increasing sequence of relative heights
\begin{equation}\label{relh}h_i:=\frac{l_{j_i}}{\langle\lambda,\beta_{j_i}^\vee\rangle}\,\in[0,1)\cap{\mathbb{Q}}\,,\;\;\;\;\;i=1,\ldots,s\,.\end{equation}
Let $0<b_2<\cdots<b_t<1$ be the distinct nonzero values in the set $\{h_1,\ldots,h_s\}$, and let $b_1:=0$, $b_{t+1}:=1$. For $k=1,\ldots,t$, let $I_k:=\{1 \leq i \leq s \, \mid \, h_i=b_k\}$.

Recall the path $\Pi(w,A)$ in $\QB(W)$ defined in~\eqref{eqn:admissible}. We divide this path into subpaths corresponding to the sets $I_k$, and record the last element in each subpath; more precisely, for $k=0,\ldots,t$, we define the sequence of Weyl group elements
\[w_k:=w\prod_{i\in I_1\cup\cdots\cup I_k}^{\longrightarrow} r_{j_i}\,,\]
where the non-commutative product is taken in the increasing order of the indices $i$, and $w_0:=w$. For $k=1,\ldots,t$, let $\sigma_k:=\lfloor w_k\rfloor\in W^J$. We can now define the forward map as
\[(w,A)\,\mapsto\,\eta:=((b_1,b_2,\ldots,b_t,b_{t+1});\,(\sigma_1,\ldots,\sigma_t))\,.\]
We will verify below that the image is in $\QLS(\lambda)$. 

The inverse map is constructed using the quantum Deodhar lift and the related {\em shellability property} of the quantum Bruhat graph~\cite{BFP}. 

\begin{prop}\label{qam2qls} The map $A\mapsto \eta$ constructed above is a bijection between ${\mathcal A}(w,\Gamma_{\rm lex}(\lambda))$ and $\QLS(\lambda)$. It maps the corresponding parameters in the following way:
\begin{equation*}\wt(w,A)=\wt(\eta),\: {\rm end}(w,A)=\iota(\eta,w),\: {\rm down}(w,A)=\xi(\eta,w),\: -{\rm height}(w,A)=\Deg_w(\eta).\end{equation*}
\end{prop}

We translate the formula in Theorem~\ref{chevdomqls} to the quantum alcove model via Prop.~\ref{qam2qls}.

\begin{thm}\label{chevdomqam} Let $\lambda$ be a dominant weight, $\Gamma_{\rm lex}(\lambda)$ the lex $\lambda$-chain, and let $x=wt_{\xi}\in \Wafp$. Then, in $K_{T \times \mathbb{C}^*}(\QQ)$, we have
\[\begin{split}&\OQGl{-w_\circ\lambda}\cdot\OQG{x}=\\[3mm]&\qquad \sum_{A\in{\mathcal A}(w,\Gamma_{\rm lex}(\lambda))}\sum_{{\bm{\chi}}\in\overline{{\rm Par}(\lambda)}}q^{-{\rm height}(w,A)-\langle\lambda,\xi\rangle-|{\bm{\chi}}|}\mathbf{e}^{{\rm wt}(w,A)}\OQG{{\rm end}(w,A)t_{\xi+{\rm down}(w,A)+\iota({\bm{\chi}})}}\,.\end{split}\]
\end{thm}

A similar Chevalley formula for an anti-dominant weight $\lambda$ was derived in \cite[Theorem~1]{nospcf}, also in terms of quantum LS paths. Using a similar procedure to the one above, we translate it to the quantum alcove model, as stated in Theorem~\ref{corantidom}. We work with  the lex $\lambda$-chain $\Gamma_{\rm lex}(\lambda)$, defined as the reverse of the lex $(-\lambda)$-chain in Proposition~\ref{speciallch}; note that the alcove path corresponding to the former (ending at $A_\circ-\lambda$) is the translation by $-\lambda$ of the alcove path corresponding to the latter (ending at $A_\circ+\lambda$). 

\begin{thm}\label{corantidom} Let $\lambda$ be an anti-dominant weight, $\Gamma_{\rm lex}(\lambda)$ the lex $\lambda$-chain, and let $x=wt_{\xi}\in \Wafp$. Then, in $K_{T \times \mathbb{C}^*}^{\prime}(\QQ) \subset K_{T \times \mathbb{C}^*}(\QQ)$, we have
\[\begin{split}&\OQGl{-w_\circ\lambda}\cdot\OQG{x}=\\[3mm]&\qquad\sum_{A\in{\mathcal A}(w,\Gamma_{\rm lex}(\lambda))}(-1)^{|A|}q^{-{\rm height}(w,A)-\langle\lambda,\xi\rangle}\mathbf{e}^{{\rm wt}(w,A)}\OQG{{\rm end}(w,A)t_{\xi+{\rm down}(w,A)}}\,.
\end{split}\]
\end{thm}

\subsection{The Chevalley formula for an arbitrary weight}

We now exhibit the Chevalley formula for an arbitrary weight $\lambda=\sum_{i\in I}\lambda_i\varpi_i$, for $q=1$. The proof is discussed in Section~\ref{secproof}. To state the formula, let $\overline{{\rm Par}(\lambda)}$ denote the set of $I$-tuples of partitions $\bm{\chi}=(\chi^{(i)})_{i\in I}$ such that $\chi^{(i)}$ is a partition of length at most $\max(\lambda_i,0)$.

\begin{thm}\label{genchev} Let $\lambda$ be an arbitrary weight, $\Gamma(\lambda)$ an arbitrary $\lambda$-chain, and let $x=wt_{\xi}\in \Wafp$. Then, in $K_{T\times{\mathbb C}^*}(\QQ)$, we have
\[\begin{split}&\OQGl{-w_\circ\lambda}\cdot\OQG{x}=\\[3mm]&\quad\!\!\!\!\sum_{A\in{\mathcal A}(w,\Gamma(\lambda))}\sum_{{\bm{\chi}}\in\overline{{\rm Par}(\lambda)}}(-1)^{n(A)}q^{-{\rm height}(w,A)-\langle\lambda,\xi\rangle-|{\bm{\chi}}|}\mathbf{e}^{{\rm wt}(w,A)}\OQG{{\rm end}(w,A)t_{\xi+{\rm down}(w,A)+\iota({\bm{\chi}})}}\,,
\end{split}\]
where $n(A)$, for $A=\{j_1<\cdots<j_s\}$, is the number of negative roots in $\{\beta_{j_1},\ldots,\beta_{j_s}\}$.
\end{thm}

\begin{exa}{\rm 
Assume that $\mathfrak{g}$ is of type $A_2$, and $\lambda=\varpi_1-\varpi_2$. 
Assume that $\mathfrak{g}$ is of type $A_2$, and $\lambda=\omega_1-\omega_2$. 
Then, $\Gamma(\lambda):=(\alpha_{1},\,-\alpha_{2})$ is a $\lambda$-chain of roots.
Assume that $w = s_{1}=s_{\alpha_{1}}$. In this case, we see that
${\mathcal A}(s_{1},\Gamma(\lambda))=\bigl\{\emptyset,\,\{1\},\,\{2\},\,\{1,2\}\bigr\}$, 
and we have the following table. 
\begin{equation*}
\begin{array}{c||c|c|c|c|c}
A & n(A) & {\rm height}(s_{1},A) & {\rm wt}(s_{1},A) & {\rm end}(s_{1},A) & {\rm down}(s_{1},A) \\ \hline\hline
\emptyset & 0 & 0 & s_{1}\lambda & s_{1} & 0 \\ \hline
\{1\} & 0 & 1 & \lambda & e & \alpha_{1}^{\vee} \\ \hline
\{2\} & 1 & 0 & s_{1}\lambda & s_{1}s_{2} & 0 \\ \hline
\{1,2\} & 1 & 1 & \lambda & s_{2} & \alpha_{1}^{\vee}
\end{array}
\end{equation*}
Also we can identify $\overline{{\rm Par}(\lambda)}$ with $\mathbb{Z}_{\ge 0}$. 
Therefore we obtain 
\begin{equation*}
\begin{split}
 & \OQGl{-w_\circ\lambda} \cdot \OQG{s_{1}t_{\xi}} = \\[3mm]
 & \hspace*{5mm} 
  \sum_{ m \in \mathbb{Z}_{\ge 0} } q^{ -\langle \lambda,\xi \rangle-m}
  \biggl\{%
  \underbrace{
  \mathbf{e}^{s_1\lambda}
  \OQG{ s_1 t_{\xi+m\alpha_{1}^{\vee}} }}_{A=\emptyset} 
  + 
  \underbrace{
  q^{-1} \mathbf{e}^{\lambda}
  \OQG{ t_{\xi+\alpha_{1}^{\vee}+m\alpha_{1}^{\vee}} }}_{A=\{1\}} \\
 & \hspace*{40mm}
   + 
   \underbrace{(-1)\mathbf{e}^{s_{1}\lambda}
   \OQG{ s_1s_2 t_{\xi+m\alpha_{1}^{\vee}} }}_{A=\{2\}}
   +
   \underbrace{ (-1)q^{-1} \mathbf{e}^{\lambda}
   \OQG{ s_2 t_{\xi+\alpha_{1}^{\vee}+m\alpha_{1}^{\vee}} }}_{A=\{1,2\}} \biggr\}. 
\end{split}
\end{equation*}}
\end{exa}

\section{The quantum $K$-theory of flag varieties}\label{qkgb}

Y.-P. Lee defined the (small) {\em quantum $K$-theory} of a smooth projective variety $X$, denoted by $QK(X)$ \cite{leeqkt}. This is a deformation of the ordinary $K$-ring of $X$, analogous to the relation between quantum cohomology and ordinary cohomology. The deformed product is defined in terms of certain generalizations of {\em Gromov-Witten invariants} (i.e., the structure constants in quantum cohomology), called {\em quantum $K$-invariants of Gromov-Witten type}. 

In order to describe the (small) $T$-equivariant quantum $K$-algebra $QK_T(G/B)$, for the finite-dimensional flag manifold $G/B$, we associate a variable $Q_i$ to each simple coroot $\alpha_i^{\vee}$, 
and let $\bZ[Q]:=\bZ[Q_1,\dots,Q_r]$. 
Given $\mu=d_1 \alpha_1^\vee+\cdots+d_r\alpha_r^\vee$ in $Q^{\vee,+}$, let $Q^\mu:=Q_1^{d_1}\cdots Q_r^{d_r}$. 
Let $\bZ[Q][P]:=\bZ[Q]\otimes_{{\mathbb Z}} \bZ[P]$, where the group algebra $\mathbb{Z}[P]$ of $P$ was defined at the beginning of Section~\ref{csi}. 
As a $\bZ[Q][P]$-module, $QK_T(G/B)$  is defined as
$K_T(G/B) \otimes_{{\mathbb Z}[P]} \bZ[Q][P]$. 
The algebra $QK_T(G/B)$ has a {$\mathbb{C}[P]$-basis given by} the classes $[{\mathcal O}^w]$ of the structure sheaves of (opposite) Schubert varieties in $G/B$, for $w \in W$. 

It is proved in \cite{kat1} (see also \cite{kat2}) that there exists a $\mathbb{C}[P]$-module isomorphism from $QK_T(G/B)$ onto $K_{T}^{\prime}(\QQ)$ that respects the quantum multiplication in $QK_T(G/B)$ and the tensor product in $K_{T}^{\prime}(\QQ)$;
in particular, it respects the quantum multiplication with a line bundle $\mathcal{O}_{G/B}(-\varpi_{i})$ and the tensor product with a line bundle $\OQGl{w_{\circ} \varpi_{i}}$, for $i \in I$.
Note that this isomorphism sends each (opposite) Schubert class $[{\mathcal O}^w] Q^{\mu}$ in $QK_T(G/B)$ to the corresponding semi-infinite Schubert class $\OQG{wt_{\mu}}$ in $K_{T}^{\prime}(\QQ)$ for $w \in W$ and $\mu \in Q^{\vee,+}$. The above result and the formula in Theorem~\ref{corantidom} imply an important conjecture in~\cite{LP}: the Chevalley formula for $QK_T(G/B)$. 

\begin{thm}\label{qkchev} Let $i \in I$, and fix a $(-\varpi_i)$-chain of roots $\Gamma(-\varpi_i)$. Then, in $QK_{T}(G/B)$, we have
\[\begin{split}&[{\mathcal O}^{s_i}]\cdot [{\mathcal O}^w] = \\[3mm]&\quad=
(1-\mathbf{e}^{w(\varpi_i) - \varpi_i})  [{\mathcal O}^w] +
\sum_{A\in{\mathcal A}(w,\Gamma(-\varpi_i))\setminus\{\emptyset\}}
(-1)^{|A|-1} \,Q^{{\rm{down}}(w,A)}\mathbf{e}^{-\varpi_i-{\rm{wt}}(w,A)} [{\mathcal O}^{{\rm end}(w,A)}].\end{split}\]
\end{thm}

Let us now turn to the type $A$ flag manifold $Fl_{r+1}=SL_{r+1}/B$ and its (non-equivariant) quantum $K$-theory $QK(Fl_{r+1})$. In~\cite{lamqgp}, the first author and Maeno defined the so-called {\em quantum Grothendieck polynomials}. According to~\cite[Theorem~6.4]{lamqgp}, whose proof is based on intricate combinatorics, the mentioned polynomials multiply precisely as stated by the above Chevalley formula. As this formula determines the multiplicative structure of $QK(Fl_{r+1})$, we derive the following result, settling the main conjecture in~\cite{lamqgp}.

\begin{thm}\label{qgroth} The quantum Grothendieck polynomials represent Schubert classes in $QK(Fl_{r+1})$. 
\end{thm}

Given a {\em degree} $d=(d_1,\ldots, d_r)$, let $N_{s_i,w}^{v,d}$ be the coefficient of $Q_1^{d_1}\cdots Q_r^{d_r}[{\mathcal O}^v]$ in the expansion of $[{\mathcal O}^{s_i}]\cdot [{\mathcal O}^w]$ in $QK(Fl_{r+1})$ for $1 \leq i \leq r$. Based on Theorem~\ref{qkchev} and results in \cite{lenfmp}, we describe more explicitly the quantum $K$-Chevalley coefficients for $Fl_{r+1}$.

\begin{thm}
For $QK(Fl_{r+1})$, we always have $N_{s_i,w}^{v,d}\in\{0,\pm 1\}$ for $1 \leq i \leq r$. For every $i,\,v$ and parabolic coset $\sigma W_{I\setminus \{i\}}$, there are unique $d$ and $w\in \sigma W_{I\setminus \{i\}}$  (they can be constructed explicitly), such that $N_{s_i,w}^{v,d}=\pm 1$ (the sign is as in Theorem~{\rm \ref{qkchev}}), where $I = \{1, \ldots, r\}$. Moreover, we determined the maximum degree in the Chevalley formula for $QK(Fl_{r+1})$. 
\end{thm}

\section{On the proof of Theorem~\ref{genchev}. Quantum Bruhat operators}\label{secproof}

We discuss the main idea underlying the proof of Theorem~\ref{genchev}. Given an arbitrary weight $\lambda=\sum_{i\in I}\lambda_i\varpi_i$, we express it as $\lambda=\lambda^++\lambda^-$, where $\lambda^+:=\sum_{\lambda_i>0}\lambda_i\varpi_i$ is dominant and $\lambda^-:=\sum_{\lambda_i<0}\lambda_i\varpi_i$ is anti-dominant. By concatenating the chains of roots $\Gamma_{\rm lex}(\lambda^+)$ and $\Gamma_{\rm lex}(\lambda^-)$, we obtain a chain $\Gamma(\lambda)$. By Theorems~\ref{chevdomqam} and \ref{corantidom}, we can write the Chevalley formula for $\lambda$ as in Theorem~\ref{genchev}, based on $\Gamma(\lambda)$. However, this is not a $\lambda$-chain, because the corresponding alcove path from $A_\circ$ to $A_{-\lambda}$ is not reduced. Thus, the problem is to deform a non-reduced alcove path to a reduced one; this also takes care of possible cancellations of terms. Moreover, as we want to prove that the formula is independent of the $\lambda$-chain, we need to relate two reduced alcove paths.

By Remark~\ref{remredword}, deforming alcove paths amounts to relating reduced words for elements of $\Waf$ via Coxeter relations. This approach was used in \cite{LP} to prove the Chevalley formula for $K_T(G/B)$ in terms of the alcove model. The idea was to express the formula as a composition of operators $R_\alpha$, for $\alpha$ ranging over the sequence of (positive or negative) roots encoding an alcove path (not necessarily reduced). In this setup, the independence of the formula from the considered alcove path amounts to the fact that the family $\{R_\alpha : \alpha\in\Phi\}$ satisfies the {\em Yang-Baxter equation}. This means that $R_{-\alpha}$ is the inverse of $R_\alpha$ and, for any pair of roots
$\alpha,\beta\in\Phi$ such that $(\alpha,\beta)\leq 0$, we have
\begin{equation}\label{ybe}
R_{\alpha} R_{s_\alpha(\beta)} 
R_{s_\alpha s_\beta(\alpha)}\cdots
R_{s_\beta(\alpha)} R_\beta
=
R_\beta R_{s_\beta(\alpha)} \cdots
R_{s_\alpha s_\beta(\alpha)} R_{s_\alpha(\beta)} R_{\alpha}\,.
\end{equation}

For simplicity, we now restrict to the non-equivariant version of Theorem~\ref{genchev}, and refer to \cite[Sections~10,~13,~14]{LP} for addressing the equivariant version. Using the notation in Section~\ref{qkgb}, the relevant operators for our proof are defined as follows. For $\alpha\in\Phi^+$, let $R_\alpha:=1+Q_\alpha$, where $Q_\alpha$ is the quantum Bruhat operator in \cite{BFP}; this is defined on the group algebra $\bZ[Q][W]$ of $W$ over $\bZ[Q]$ by:
$$
Q_{\alpha}(w)=\left\{ \begin{array}{ll}  w s_\alpha &\mbox{if}\;\; w\lessdot ws_\alpha, \\ 
Q^{\alpha^\vee}\,w s_\alpha &\mbox{if}\;\; w \triangleleft ws_\alpha,\\ 
0 &\mbox{otherwise\,.} \end{array}
\right.  
$$
By \cite[Corollary~4.4]{BFP}, \cite[Lemma~9.2]{LP}, $\{R_\alpha:\alpha\in\Phi\}$ satisfies the Yang-Baxter equation.

However, the proof of Theorem~\ref{genchev} is not yet concluded. This is because, in order to factor our Chevalley formula, we need to replace $R_{-\alpha}=(R_\alpha)^{-1}=1-Q_\alpha+Q_\alpha^2-\ldots$ with $R_{-\alpha}':=1-Q_\alpha$, for $\alpha\in\Phi^+$; indeed, unlike in the setup of \cite{LP}, $Q_\alpha^2$ is not always $0$. In order to address this problem, we need a finer analysis of the Yang-Baxter property~\eqref{ybe}. For $\lambda$-chains with dominant $\lambda$, this was done in \cite{lalurc}, via certain combinatorial moves (the {\em quantum Yang-Baxter moves}); their role is to biject the paths in $\QB(W)$ indexing the terms on the two sides of~\eqref{ybe}. The proof of Theorem~\ref{genchev} is concluded using a generalization of this work to arbitrary alcove paths.

\end{document}